\documentclass[12pt]{article}
\usepackage{amsmath}
\usepackage{amsthm}
\usepackage{amsfonts}
\usepackage{amssymb}
\usepackage{amstext}
\usepackage{amsopn}
\input epsf

\newtheorem{thm}{Theorem}

\newtheorem{lemma}{Lemma}

\begin{document}
 \def\today{20~November~2009, revised 21 January 2010} 
\title{\bf 
Meromorphic solutions of a third order nonlinear differential equation}
\author{\textsc{Robert Conte${}^{1,2}$\thanks{ 
Partially supported by PROCORE - France/Hong Kong joint research grant
F-HK29/05T and RGC grant HKU 703807P.
\mbox{\hspace{4truecm}}} 
and Tuen-Wai Ng \epsfxsize 15mm \epsffile{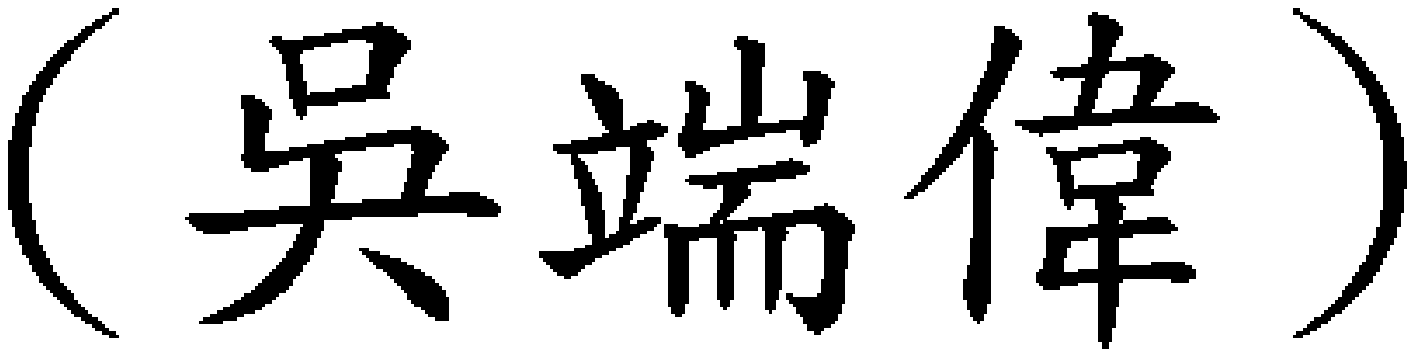} ${}^{1}$\thanks{
Partially supported by PROCORE - France/Hong Kong joint research grant 
F-HK29/05T and RGC grant HKU 703807P.\mbox{\hspace{0.9truecm}}}}}

\date{}
\maketitle

\centerline{\textbf{\today}}

\begin{figure}[b]
\rule[-2.5truemm]{5cm}{0.1truemm}\\[2mm]
{\footnotesize  
2000{\it Mathematics Subject Classification: Primary} 30D35.
\par\noindent 1.
Department of Mathematics,
The University of Hong Kong,
Pokfulam Road.
\smallskip
\par\noindent 2.
LRC MESO, 
Centre de math\'ematiques et de leurs applications (UMR 8536) et CEA-DAM, 
\\ \'Ecole normale sup\'erieure de Cachan, 61, avenue du Pr\'esident Wilson,
\\ F--94235 Cachan Cedex, France.
\smallskip
\par\noindent E-mail: Robert.Conte@cea.fr, ntw@maths.hku.hk 
\par
}

\end{figure}

\begin{quotation}
\noindent{\sc\textbf{Abstract}.  
We prove that all the
meromorphic solutions of the nonlinear differential equation
$c_0u''' + 6u^4 + c_1u'' + c_2u u' + c_3u^3 + c_4u'+ c_5u^2 + c_6u +c_7=0$ 
are elliptic or degenerate elliptic,
and we build them explicitly.}
\end{quotation}

\def \res {\mathop{\rm res}\nolimits}


\vfill\eject

\section{Introduction}
\label{sectionIntro}

When a system is governed by an autonomous nonlinear algebraic 
partial differential equation (PDE),
it frequently admits some permanent profile structures
such as fronts, pulses, kinks, etc
\cite{vS2003},
and usually these profiles are mathematically some 
single-valued solutions 
of the travelling wave reduction $(x,t) \to x-ct$
of the PDE to an ordinary differential equation (ODE).
The physical motivation of the present work
is to find such solutions in closed form.
Since this is a difficult mathematical problem,
we restrict here to a simple case (a third order nonlinear ODE)
and solve it completely. 
The method we used here is a refinement of Eremenko's method 
used in \cite{Eremenko1982} as well as \cite{EremenkoKS} 
and \cite{ELN} which is based on the local singularity analysis 
of the meromorphic solutions of the given differential equations 
as well as the zero distribution and growth rate of the meromorphic 
solutions by using Nevanlinna theory. 
This is a very powerful method. 
For example, it was used by Eremenko \cite{EremenkoKS} 
to characterize all meromorphic traveling wave solutions of the 
Kuramoto-Sivashinsky (KS) equations. 
In fact, Eremenko showed that all the meromorphic traveling wave solutions
 of the KS equations belong to the class $W$ (like Weierstrass),
which consists of 
elliptic functions and their successive degeneracies,
i.e.: 
elliptic functions, 
rational functions of one exponential $\exp(kz), k\in\mathbb{C}$
and rational functions of $z$.

In general, even if we know that the solutions belong to the class $W$, 
it is still difficult to find their explicit form. 
To overcome this problem, 
we shall apply the subequation method introduced in \cite{MC2003}
and developed in \cite{CM2009}.
In order to emphasize the method, 
we will choose a test equation according to the following criteria:
\begin{enumerate}
\item
to have a small differential order $n$,

\item
to have only nonrational Fuchs indices,
apart from the ever present $-1$ index,

\item
to be of the form
$u^{(n)}=P(u^{(n-1)},\dots,u',u)$,
with $P$ a polynomial of its arguments,

\item
to have movable poles of order one,

\item
to be \textit{complete} in the classical sense \cite{PaiBSMF}
(see details in \cite[p.~122]{Cargese1996Conte})
i.e.~to include all admissible nondominant terms,

\end{enumerate}

The requirement for nonrational Fuchs indices sets $n \ge 3$.
Let us take the complete autonomous
third order polynomial ODE with simple poles, 
\begin{eqnarray}
& &
d_0 u''' + d_1 u u'' + d_2 {u'}^2 + d_3 u^2 u' + d_4 u^4
\nonumber \\ & &
 +c_1 u'' +c_2 u u' +c_3 u^3 
 +c_4 u'+c_5 u^2
 +c_6 u
 +c_7=0.
\label{eqODE3Complete}
\end{eqnarray}
This equation is indeed complete in the sense that it
includes all polynomial terms having a singularity degree at most
equal to four,
as seen from the generating function
\begin{eqnarray}
& &
\frac{1}{(1-t u)(1 -t^2 u')(1-t^3 u'')}
\nonumber\\ & &
=1 + u t
+(u^2+u') t^2 + (u^3 + u u' + u'') t^3 
\nonumber\\ & & \phantom{1234567890}
+(u^4 + u^2 u' +{u'}^2 + u u'') t^4+O(t^5).
\end{eqnarray}
Let us choose one particular set of dominant terms 
(the ones with coefficients $d_j$, which have quadruple poles)
so as to enforce from the beginning
the condition that the Fuchs indices be nonrational. 
After setting $c_3=0$ without loss of generality,
our test equation will be normalized as
\begin{eqnarray}
& &
c_0 u''' + 6 u^4
 +c_1 u'' +c_2 u u' 
 +c_4 u'+c_5 u^2
 +c_6 u
 +c_7=0,
\label{eqODE3}
\end{eqnarray}

Let $u$ be a meromorphic solution of the ODE (\ref{eqODE3}).
We first check that if $u$ has a movable 
pole at $z=z_0$,
then $u$ has only three distinct Laurent series expansions at $z_0$. 
Note that if $z_0$ is a pole of $u$, it must be a simple pole. 
Therefore, in a neighbourhood of $z=z_0$, 
the Laurent series of the meromorphic solution $u$ is of the form 
\begin{eqnarray}
& &
u(z)= u_{-1}(z-z_0)^{-1} +u_0 +u_1 (z-z_0)+\cdots,\ u_{-1}\not=0.
\label{eqLaurent}
\label{eqODE3Laurent}
\end{eqnarray}
Denote $a$ any one 
of the cubic roots of $c_0$. 
Substituting the above Laurent series into the ODE (\ref{eqODE3})
and balancing the leading terms, 
we obtain $u_{-1}=a$,
and $u_0=(-2 c_1 a + c_2 a^2)/(24 c_0)$. 
We are going to prove that there are at most 
three distinct Laurent series expansions at $z_0$.
If one linearizes the ODE (\ref{eqODE3}) 
around the movable singularity $z=z_0$ \cite[p.~114]{Cargese1996Conte},
the resulting linear ODE has the Fuchsian type at $z_0$,
and its three Fuchs indices $r$ are defined by
\begin{eqnarray}
& &
(r+1) (r^2-7r+18)=0.
\end{eqnarray}

Hence, the Fuchs indices
are equal to $r=-1,(7 \pm \sqrt{-23})/2$.
Because of the absence of any positive integer 
in the set of values of $r$,
all other cofficients $u_i$ are uniquely determined 
\cite[p.~90]{Cargese1996Conte}
by the leading coefficient $u_{-1}$. 
Hence, there are at most three meromorphic functions     
with poles at $z=z_0$ satisfying the ODE (\ref{eqODE3}).

We shall study the third order nonlinear 
differential equation (\ref{eqODE3}) 
and show that all meromorphic solutions of this differential 
equation belong to the class $W$. 
More specifically, our results are the following.

\begin{thm} 
If the ODE (\ref{eqODE3}) has a particular meromorphic solution $u$, 
then $u$ belongs to the class $W$.  
Moreover, a necessary and sufficient condition for the ODE (\ref{eqODE3}) 
to admit a particular meromorphic solution is
to belong to the following list,
\begin{eqnarray}
& &
S_{3a}:\
c_1,c_6=\mbox{\rm arbitrary},\
c_2=0,\ c_5=0,\ c_7=0,\ c_4=\frac{c_1^2}{12 c_0};\
\\
& &
S_{3b}:\
c_5,c_6=\hbox{\rm arbitrary},\
c_1=0,\ c_2=0,\ c_4=0,\ c_7=\frac{c_5^2}{128};\
\\
& &
S_{2A}:\
c_1,c_4=\hbox{\rm arbitrary},\
c_2=0,\ 
c_5=\frac{c_1^2 - 12 a^3 c_4}{4 a^4},\
c_6=-\frac{c_1 (c_1^2 + 36 a^3 c_4)}{144 a^6},\
\nonumber \\ & & \phantom{12345}
c_7=\frac{(12 a^3 c_4 -c_1^2)(36 a^3 c_4-11 c_1^2)}{1536 a^8};\
\\
& &
S_{2B}:\
c_1,c_2=\hbox{\rm arbitrary},\
c_4= \frac{44 c_1^2 + 8 a c_1 c_2 -   a^2 c_2^2}{144 a^3},\
\nonumber \\ & & \phantom{12345}
c_5=\frac{-32 c_1^2 -24 a c_1 c_2 - 7 a^2 c_2^2}{48 a^4},\
c_6=-\frac{(c_1+a c_2)(12 c_1^2 + 6 a c_1 c_2 + a^2 c_2^2)}{144 a^6},\
\nonumber \\ & & \phantom{12345}
c_7=-\frac{(4 c_1 + 3 a c_2)(48 c_1^2 +20 a c_1 c_2 + a^2 c_2^2)}{55296 a^7};
\\
& &
S_{1}:\
c_1,c_2,c_4,c_5=\hbox{\rm arbitrary},\
\nonumber \\ & & \phantom{1234}
1152 a^6 c_6=-56 c_1^3 + 60 a c_1^2 c_2 -18 a^2 c_1 c_2^2
             +a^3 c_2^3 +288 a^3 c_1 c_4 
\nonumber \\ & & \phantom{1234567890123}
             -144 a^4 c_2 c_4-96 a^4 c_1 c_5 +48 a^5 c_2 c_5,
\nonumber \\ & & \phantom{1234}
2^{13} 3^2 a^8 c_7=-176 c_1^4 +128 a c_1^3 c_2 +24 a^2 c_1^2 c_2^2 
                   -32 a^3 c_1 c_2^3
\nonumber \\ & & \phantom{1234567890123}
                   +5 a^4 c_2^4 +2688 a^3 c_1^2 c_4
                   -1536 a^4 c_1 c_2 c_4 +96 a^5 c_2^2 c_4
                   -6912 a^6 c_4^2 
\nonumber \\ & & \phantom{1234567890123}
                   +128 a^4 c_1^2 c_5 -512 a^5  c_1 c_2 c_5 
                   + 224 a^6  c_2^2  c_5 
\nonumber \\ & & \phantom{1234567890123}
                   + 4608 a^7  c_4 c_5 + 2304 a^8  c_5^2.
\end{eqnarray}

\end{thm}

\bigskip

We shall apply Eremenko's method \cite{EremenkoKS} 
to prove the first part of Theorem 1.
Here, we shall assume the readers are familiar with the standard terminology 
and results of Nevanlinna theory. 
The standard reference of this theory are \cite{Hayman64} 
and \cite{Laine-book,Nevanlinna-book} 
(see also \cite{EremenkoKS} for a quick introduction). 
Our argument is slightly different from that of Eremenko 
and it makes use of the following version of Clunie's Lemma 
(\cite[Lemma 2.4.2]{Laine-book}, see also \cite{YZ}).

\begin{lemma}  Let $f$ be a transcendental meromorphic solution of 
$$f^nP(z,f)=Q(z,f),$$
where $P(z,f)$ and $Q(z,f)$ are polynomials in $f$ 
and its derivatives with meromorphic coefficients 
$\{a_\lambda|\lambda \in I\}$ such that 
$m(r,a_{\lambda})=S(r,f)$ for all $\lambda \in I$. 
If the total degree of $Q(z,f)$ as a polynomial in $f$ 
and its derivatives is less than or equal to $n$, then 
$$m(r,P(r,f))=S(r,f).$$
\label{lemma2}
\end{lemma}

Now let $u$ be a function meromorphic in the complex plane 
which satisfies the above ODE (\ref{eqODE3}). 
If $u$ is rational, then we are done.
So suppose $u$ is a transcendental meromorphic solution 
of equation (\ref{eqODE3}), then we have 
\begin{eqnarray}
& &
-6 u^4=c_0 u''' 
 +c_1 u'' +c_2 u u'
 +c_4 u'+c_5 u^2
 +c_6 u
 +c_7,
\end{eqnarray}
Take $f=u, P=u, n=3$ and apply Clunie's lemma (Lemma 1) 
to the above equation, 
we conclude that $m(r,u)=S(r,u)$ and hence $(1-o(1))T(r,u)=N(r,u)$. 
We claim that $u$ must have infinitely many poles. 
Assume it is not the case, then $N(r,u)= O(\log r)$.
Therefore, $T(r,u)=O(\log r)$ which is impossible as $u$ is transcendental. 

Secondly, we  prove that if $u$ is a transcendental 
meromorphic solution, 
then $u$ is a periodic function. 
Recall that there are at most three meromorphic functions 
with poles at $z=z_0$ satisfying the ODE (\ref{eqODE3}). 
Now let $z_j, j=1,2,3,\cdots$ be the poles of $u(z)$, 
then the functions $w_j(z) = u(z+z_j-z_0)$ are meromorphic 
solutions of the ODE (\ref{eqODE3}) with a pole at $z_0$. 
Thus some of them must be equal. 
Consequently, $u$ is a periodic function. 

Without loss generality, we may assume that $u$ has a period of $2\pi i$. 
Let
$D=\{z:0\le {\rm Im} z <2\pi\}$. If $u$ has more than three poles in $D$, 
then by the previous argument, 
we can conclude that $u$ is periodic in $D$ 
and hence it is indeed an elliptic function and we are done.

Now suppose $u$ has at most three poles in $D$.
Since $u$ is a periodic function with period $2\pi i$, we have
$N(r,u)= O(r),$ as $r\to\infty$. 
It follows from  $(1-o(1))T(r,u)=N(r,u)$ that $T(r,u)=O(r)$. 
By Nevanlinna's First Fundamental Theorem, we know that 
for any $a\in\mathbb{C}$, $N(r,1/(u -a)) =O(r)$ as $r\to\infty$. 
By the periodicity of $u$, we conclude that $u$ take each $a$ finitely many
times in $D$. Hence, the function $R(z)=u(\ln z)$ is a
single-valued analytic function in the punctured plane
$\mathbb{C}\backslash \{0\}$ and takes each $a \in \mathbb{C}$ finitely many times. 
It follows that $0$ is a removable singularity of $R$ and $R$ 
can then be extended to a meromorphic function on $\mathbb{C}$. Hence $R$ is a rational function 
as it takes each complex number finitely many times. 
Therefore, $u(z)=R(e^z)$ belongs to the class $W$ and this completes the proof of 
the first part of Theorem 1.

\medskip
\noindent
{\bf Remark.} 
{}From the above proof, we notice that if $u$ is an elliptic solution, 
then $u$ has at most three (simple) poles in each fundamental polygon 
$\Omega$. 
Recall that the residue of $u$ at any pole must be one of 
$a,\omega a,\omega^2 a$ where $\omega$ is the cubic root of unity. 
Since the sum of the residues of all the poles in any fundamental 
polygon $\Omega$ is zero, 
$u$ must have three distint simple poles in $\Omega$ 
and hence we have three distinct Laurent series at $z_0$.

\noindent
{\bf Remark.} 
If $u(z)=R(e^{kz})$ where $R$ is some rational function, 
then $R$ has at most three (simple) poles in $\mathbb{C}\backslash \{0\}$. 
We are going to show that $R$ cannot have a pole at $0$. 
Suppose we write 
$u(z)=R(Z)= r_0/Z^n + \sum_{i=1}^{3} r_i/(Z-Z_i) + q(Z)$, 
where $q$ is a polynomial in $Z= e^{kz}$. 
Substituting $u(z)=R(Z)$ into ODE (\ref{eqODE3}) and letting $Z$ tend to infinity, 
we can conclude that $q$ equals to some constant $C$. 
Now letting $Z$ tend to $0$, we can deduce that $r_0=0$.
Hence, $u(z)= \sum_{i=1}^{3} r_i/(e^{kz} -Z_i) + C$, 
where $Z_i, C \in \mathbb{C}$.  
Finally, if $u$ is rational, 
then $u$ will have at most three (simple) poles in $\mathbb{C}$.  
Similarly, we can show that $u$ must be of the form $\sum_{i=1}^{3}r_i/(z-z_i) + C$, 
where $r_i, C \in \mathbb{C}$.

\section{Explicit solutions in the class $W$}
\label{section2}

Let us determine the constraints on the coefficients 
$c_j$ of (\ref{eqODE3}) for meromorphic solutions to exist,
and let us determine all these meromorphic solutions
in closed form.
According to section \ref{sectionIntro},
these solutions are necessarily elliptic or degenerate of elliptic
(i.e.~rational in one exponential $e^{k z}, k \in \mathbb{C}$
or rational in $z$),
i.e.~they belong to the class $W$.

If the meromorphic solution is elliptic,
by a classical theorem,
the sum of the residues of the three Laurent series
for $u$, Eq.~(\ref{eqLaurent}),
must vanish,
and similarly for any rational function of $u,u',u''$.
These necessary conditions \cite{Hone2005} 
are first established in section
\ref{sectionResidueConditions}.

If the solution is elliptic,
one knows the elliptic orders of $u$ and $u'$,
they are respectively equal to 
three (three simple poles) 
and six (three double poles).
Therefore,
by a classical theorem of Briot and Bouquet
\cite[p.~277]{BriotBouquet},
\cite[part II, chap.~IX p.~329]{HalphenTraite},
\cite[p.~424]{Hille}
the elliptic solution obeys a first order algebraic equation
whose degree in $u'$ is the order of $u $ (three)
and   degree in $u $ is the order of $u'$ (six),
\begin{eqnarray}
& &
F(u,u') \equiv
 \sum_{k=0}^{m} \sum_{j=0}^{2m-2k} a_{j,k} u^j {u'}^k=0,\ 
a_{0,m}\not=0,\ 
\label{eqsubeqODEOrderOnePP}
\end{eqnarray}
with $m=3$.
The complex constants $a_{j,k}$, with $a_{0,m}\not=0$, 
are then determined by the algorithm presented in 
\cite{MC2003},
i.e.~by requiring each of the three Laurent series
(\ref{eqLaurent}) to obey (\ref{eqsubeqODEOrderOnePP}).
The search for all third degree subequations 
(\ref{eqsubeqODEOrderOnePP})
obeyed by the three Laurent series
is performed in section \ref{sectionsubeqdeg3}.

As to those solutions of (\ref{eqODE3}) which are
degenerate of elliptic,
they also obey a first order equation (\ref{eqsubeqODEOrderOnePP}),
whose degree $m$ is at most three.
Because of the singularity structure of (\ref{eqODE3})
(three \textit{distinct} Laurent series),
any $m$-th degree subequation, $1 \le m \le 3$,
must have $m$ \textit{distinct} Laurent series.
The search for all second or first degree subequations 
(\ref{eqsubeqODEOrderOnePP})
is performed in sections (\ref{sectionsubeqdeg2})
and (\ref{sectionsubeqdeg1}).

Let us first establish all these first order subequations.
Their general solution may be either singlevalued 
(and hence in class $W$) or multivalued.
The explicit integration of the singlevalued subset
will provide as a final output
all the meromorphic solutions of (\ref{eqODE3})
in closed form.

\subsection{Residue conditions}
\label{sectionResidueConditions}

If (\ref{eqODE3}) admits an elliptic solution,
it is necessary that,
for any rational function of $u$ and its derivatives,
the sum of the residues inside a period parallelogram be zero,
\begin{eqnarray}
& &
\forall k \in \mathbb{N}\
\forall n \in \mathbb{N}\:
\res \sum_{i=1}^{3} \left(u^{(k)}\right)^n=0.
\label{eqresidueconditions}
\end{eqnarray}
The first conditions are
\begin{eqnarray}
& &
\left\lbrace
\begin{array}{ll}
\displaystyle{
k=0,n=2:\ c_2=0,
}\\ \displaystyle{
k=0,n=3:\ c_4 =\frac{c_1^3}{12 a^3},
}\\ \displaystyle{
k=0,n=5:\ c_1 c_5 =0,
}\\ \displaystyle{
k=0,n=7:\ c_1 c_7=0,
}\\ \displaystyle{
k=1,n=4:\ (c_6 (c_5^2-128 c_7)=0 \hbox{ if } c_1=0),\
     (c_7 (c_1^3+36 a_0^2 c_6)=0 \hbox{ if } c_1\not=0).
}
\end{array}
\right.
\end{eqnarray}

When the computation is limited to $k \le 4,n \le 10$,
this defines three and only three distinct 
sets of fixed coefficients for a possible elliptic solution,
\begin{eqnarray}
& &
c_2=0,\ c_1=0,\ c_4=0,\ c_6 \not=0,\ c_7=\frac{c_5^2}{128},\ 
\label{eqrell1K0}
\label{eqrella}
\\
& &
c_2=0,\ c_1=0,\ c_4=0,\ c_6=0,\                              
\label{eqrell1K1}
\label{eqrellb}
\\
& &
c_2=0,\ c_1\not=0,\ c_4=\frac{c_1^2}{12 a^3},\ c_5=0,\ c_7=0.
\label{eqrell2K3}                                            
\label{eqrellc}                                            
\end{eqnarray}

\subsection{Subequations of degree three}
\label{sectionsubeqdeg3}

Denoting $\omega_k,k=1,2,3$, cubic roots of unity,
each such subequation has the necessary form
\begin{eqnarray}
& &
F(u,u') \equiv
-(\omega_1 a u' + u^2)(\omega_2 a u' + u^2)(\omega_3 a u' + u^2)
\nonumber \\ & & \phantom{12345678}
+b_1 {u'}^2 u + b_2 u' u^3 + b_3 u^5 
+b_4 {u'}^2   + b_5 u' u^2 + b_6 u^4 
\nonumber \\ & & \phantom{12345678}
              + b_7 u' u   + b_8 u^3 
              + b_9 u'     + b_b u^2 
                           + b_a u   
                           + b_c     
                           + b_0=0,
\label{eqsubeqODEOrderOnePP3Repeat}
\end{eqnarray}
with all $\omega_j$ distinct 
and the additional condition to be irreducible.

The first order third degree subequation is precisely defined as 
\begin{eqnarray}
& &
F(u,u') \equiv
-a^3 {u'}^3                -     u^6 
+b_1 {u'}^2 u + b_2 u' u^3 + b_3 u^5 
+b_4 {u'}^2   + b_5 u' u^2 + b_6 u^4 
\nonumber \\ & & \phantom{12345678}
              + b_7 u' u   + b_8 u^3 
              + b_9 u'     + b_b u^2 
                           + b_a u   
                           + b_c     
                           + b_0=0.
\label{eqsubeqODEOrderOnePP3}
\end{eqnarray}
The algorithm \cite{MC2003} to compute the coefficients $b_k$
is to substitute $u$ by one of
the Laurent series (\ref{eqLaurent}),
which makes
the right hand side of (\ref{eqsubeqODEOrderOnePP3}) become a Laurent series
\begin{eqnarray}
& &
F(u,u') \equiv
\sum_{j=0}^{+\infty} F_j (z-z_0)^{j-6},
\label{eqLaurentF}
\end{eqnarray}
then to solve the infinite set of equations
\begin{eqnarray}
& &
\forall a\ \forall j:\ F_j=0.
\label{eqSystemFj}
\end{eqnarray}
The practical resolution is as follows.
First, the 21 equations $F_j=0, j=0,...,6$
define a linear system for the $b_k$,
which admits a unique solution
and generates six nonlinear constraints among the six $c_k$.
By considering slightly more equations in (\ref{eqSystemFj})
(in the present case, going to $j=8$ is enough),
the set of nonlinear constraints among the $c_k$'s
admits exactly two solutions,
and all the remaining equations $F_j=0$ identically vanish.
These two solutions are
\begin{eqnarray}
& &
\left\lbrace
\begin{array}{ll}
\displaystyle{
S_{3a}:\
c_1,c_6=\hbox{arbitrary},\
c_2=0,\ c_5=0,\ c_7=0,\ c_4=\frac{c_1^2}{12 a_0},\
\label{eqgenusonecase1}
}\\ \displaystyle{
S_{3b}:\
c_5,c_6=\hbox{arbitrary},\
c_1=0,\ c_2=0,\ c_4=0,\ c_7=\frac{c_5^2}{128},
\label{eqgenusonecase2}
}
\end{array}
\right.
\end{eqnarray}
and they are identical to the two residue conditions
(\ref{eqrell2K3}) and (\ref{eqrell1K0}).

The corresponding subequations have genus one
\begin{eqnarray}
& &
(a u' + 4 k_1 u)^2 (a u' - 2 k_1 u) + (u^3 + 20 k_1^3 + k_6)^2=0,\
c_1=12 a^2 k_1,\
c_6=4 k_6,\
\label{eqsubeq3a}
\\
& &
(a u')^3 + (u^3 -3 k_5^2 u + k_6)^2=0,\
c_5=-16 k_5^2,\ 
c_6=4 k_6.
\label{eqsubeq3b}
\end{eqnarray}
The method to integrate them \cite[\S 249 p.~393]{BriotBouquet}
is to build a birational transformation to the canonical equation
of Weierstrass
\begin{eqnarray}
& &
{\wp'}^2=4 (\wp-e_1)(\wp-e_2)(\wp-e_3)=4 \wp^3 - g_2 \wp - g_3.
\end{eqnarray}
To do that, it proves convenient to introduce one of the roots $e_0$ 
of the cubic polynomial of $u(x)$ appearing as a square
in (\ref{eqsubeq3a}) and (\ref{eqsubeq3b}),
i.e.~to redefine $k_6$ by the respective relations
\begin{eqnarray}
& &
e_0^3 + 20 k_1^3 + k_6=0
\hbox{ and }
e_0^3 - 3 k_5^2 e_0 + k_6=0.
\end{eqnarray}

The subequation (\ref{eqsubeq3b})
is one of the five first order binomial equations
of Briot and Bouquet \cite[p.~122]{Cargese1996Conte},
its general solution is classical 
\begin{eqnarray}
& &
\frac{1}{u-e_0}=\frac{\wp'(z-z_0,g_2,g_3)-A}{N_1},\ 
g_2=0,\
g_3=\frac{ (e_0^2-k_5^2)^2 (e_0^2-4 k_5^2)}{243 a^6},\
\nonumber \\ & &
N_1=\frac{2 (e_0^2-k_5^2)^2}{3 a^3},\
A=\frac{e_0 (e_0^2-k_5^2)}{3 a^3}.
\label{eqBinomial3}
\end{eqnarray}

The subequation (\ref{eqsubeq3a})
has been integrated by Briot and Bouquet \cite[\S 250 p.~395]{BriotBouquet}
by introducing a function $w$ defined by
\begin{eqnarray}
& &
a u' + 4 k_1 u = \frac{u^3 - e_0^3}{u-e_0} w,
\end{eqnarray}
then by establishing the birational tranformation
\begin{eqnarray}
& &
w=\frac{a u' + 4 k_1 u}{u^2 + e_0 u + e_0^2},\
u=\frac{-3 a w w' - e_0 w^3 + 6 k_1 w^2 + 2 e_0}{2 (w^3 +1)},
\end{eqnarray}
finally by integrating the ODE for $w$,
\begin{eqnarray}
& &
w=\frac{2 k_1}{e_0}+\frac{A}{\wp - B},\
g_2=\frac{4 k_1 (k_1^3 - e_0^3)}{3 a^4},\
g_3=\frac{e_0^6 - 20 e_0^3 k_1^3 - 8 k_1^6}{17 a^6},\
\nonumber \\ & & 
g_2^3-27 g_3^2=-\frac{(8 k_1^3 + e_0^3)^3 e_0^3}{27 a^{12}}, 
A=-\frac{e_0^3 + 8 k_1^3}{3 a^2},\
B=-\frac{k_1^2}{  a^2}.
\end{eqnarray}

More generally,
birational transformations from $(u,u')$ to $(\wp,\wp')$
are obtained with an algorithm due to Poincar\'e,
implemented for instance by the command \verb+Weierstrassform+
of the computer algebra package \verb+algcurves+ \cite{MapleAlgcurves}.

\subsection{Subequations of degree two}
\label{sectionsubeqdeg2}

Let us 
define the second degree subequation as
\begin{eqnarray}
& &
F(u,u') \equiv
a^2 {u'}^2 - a u^2 u' + u^4
+b_4 u' u + b_3 u^3
+b_5 u'   + b_2 u^2
          + b_1 u
          +b_0=0,
\label{eqsubeqODEOrderTwoPP1}
\end{eqnarray}
with the additional condition to be irreducible.
Computations similar to those mentioned in section 
\ref{sectionsubeqdeg3}
provide two solutions,
\begin{eqnarray}
& &
\left\lbrace
\begin{array}{ll}
\displaystyle{
S_{2A}:\
c_1,c_4=\hbox{arbitrary},\
c_2=0,\ 
c_5=\frac{c_1^2 - 12 a^3 c_4}{4 a^4},\
}\\ \displaystyle{\phantom{12345}
c_6=-\frac{c_1 (c_1^2 + 36 a^3 c_4)}{144 a^6},\
}\\ \displaystyle{\phantom{12345}
c_7=\frac{(12 a^3 c_4 -c_1^2)(36 a^3 c_4-11 c_1^2)}{1536 a^8},\
}\\ \displaystyle{\phantom{12345}
u=v-\frac{k_1}{2},\ 
c_1=-3 a^2 k_1,\
c_4=2 a b^2 +\frac{3}{4} a k_1^2,\ 
}\\ \displaystyle{\phantom{12345}
\left(a v'-\frac{v^2-b^2}{2}\right)^2+\frac{3}{4} (v+b) (v-b) (v-k_1)^2=0,\
b \not=0,
}
\end{array}
\right.
\label{eqsubeq2sol6}
\end{eqnarray}
and
\begin{eqnarray}
& &
\left\lbrace
\begin{array}{ll}
\displaystyle{
S_{2B}:\
c_1,c_2=\hbox{arbitrary},\
c_4= \frac{44 c_1^2 + 8 a c_1 c_2 -   a^2 c_2^2}{144 a^3},\
}\\ \displaystyle{\phantom{12345}
c_5=\frac{-32 c_1^2 -24 a c_1 c_2 - 7 a^2 c_2^2}{48 a^4},\
}\\ \displaystyle{\phantom{12345}
c_6=-\frac{(c_1+a c_2)(12 c_1^2 + 6 a c_1 c_2 + a^2 c_2^2)}{144 a^6},\
}\\ \displaystyle{\phantom{12345}
c_7=-\frac{(4 c_1 + 3 a c_2)(48 c_1^2 +20 a c_1 c_2 + a^2 c_2^2)}{55296 a^7},
}\\ \displaystyle{\phantom{12345}
u=v+\frac{b}{4} + \frac{c_1}{12 a^2},\ c_2=-2 \frac{c_1}{a} + 6 a b,\
}\\ \displaystyle{\phantom{12345}
\left(a v'-\frac{v^2-b^2}{2}\right)^2+\frac{3}{4} (v+b)^3 (v-b)=0,\
b \not=0.
}
\end{array}
\right.
\label{eqsubeq2sol2}
\end{eqnarray}
For $k_1^2 \not= b^2$, 
the point transformation
\begin{eqnarray}
& &
v=k_1+\frac{1}{w},\
 w=-\frac{1}{k_1+b}-\frac{1}{k_1-b}
  +N \left(\lambda-\frac{1}{\lambda}\right),
  N^2=- \frac{b^2}{(k_1^2-b^2)^2}, 
\end{eqnarray}
maps the ODE (\ref{eqsubeq2sol6}) to the Riccati ODE
\begin{eqnarray}
& &
a N \lambda' - M \lambda-\frac{b^2}{4 (k_1^2-b^2)} (\lambda^2+1)=0,\
M^2=\frac{3 b^2}{4(k_1^2-b^2)},
\end{eqnarray}
whose general solution is a M\"obius function of one exponential so that 
$v$ is a rational function of one exponential.

For $k_1^2=b^2$, i.e.~for instance for $k_1=-b$, 
the ODE (\ref{eqsubeq2sol2}) integrates as
\begin{eqnarray}
& &
v=-b+\frac{2 b}{w},\
 w=1 + 3 \left(1 + e^{b (z-z_0)/(2 a)} \right)^2.
\end{eqnarray}

\subsection{Subequations of degree one}
\label{sectionsubeqdeg1}

These first degree subequations
\begin{eqnarray}
& &
F(u,u') \equiv
a u' + u^2 + b_1 u + b_0=0,
\label{eqsubeqODEOrderOnePP1}
\end{eqnarray}
are determined by requiring their vanishing 
when $u$ is the Laurent series (\ref{eqODE3Laurent}).
This results in
\begin{eqnarray}
& &
\left\lbrace
\begin{array}{ll}
\displaystyle{
S_{1}:\
c_1,c_2,c_4,c_5=\hbox{arbitrary},\
}\\ \displaystyle{\phantom{1234}
b_1=\frac{2 c_1 - a c_2}{12 a^2},
}\\ \displaystyle{\phantom{1234}
b_0=\frac
{44 c_1^2 - 32 a c_1 c_2 + 5 a^2 c_2^2 -144 a^3 c_4 +144 a^4 c_5}
{1152 a^4},
}\\ \displaystyle{\phantom{1234}
1152 a^6 c_6=-56 c_1^3 + 60 a c_1^2 c_2 -18 a^2 c_1 c_2^2
+a^3 c_2^3 +288 a^3 c_1 c_4 
}\\ \displaystyle{\phantom{1234567890123}
-144 a^4 c_2 c_4-96 a^4 c_1 c_5 +48 a^5 c_2 c_5,
}\\ \displaystyle{\phantom{1234}
2^{13} 3^2 a^8 c_7=-176 c_1^4 +128 a c_1^3 c_2 +24 a^2 c_1^2 c_2^2 -32 a^3 c_1 c_2^3
}\\ \displaystyle{\phantom{1234567890123}
+5 a^4 c_2^4 +2688 a^3 c_1^2 c_4
-1536 a^4 c_1 c_2 c_4 +96 a^5 c_2^2 c_4
}\\ \displaystyle{\phantom{1234567890123}
-6912 a^6 c_4^2 +128 a^4 c_1^2 c_5 -512 a^5  c_1 c_2 c_5 + 224 a^6  c_2^2  c_5 
}\\ \displaystyle{\phantom{1234567890123}
+ 4608 a^7  c_4 c_5 + 2304 a^8  c_5^2.
}
\end{array}
\right.
\label{eqSolutionOneExp}
\end{eqnarray}

The solution of this Riccati equation is 
either a rational function of one exponential or a rational function, 
\begin{eqnarray}
& &
u=
\left\lbrace
\begin{array}{ll}
\displaystyle{
-\frac{b_1}{2} + a \frac{k}{2} \coth \frac{k}{2} (z-z_0),\
k^2=\frac{b_1^2 - 4 b_0}{2 a^2}\not=0,
}\\ \displaystyle{
-\frac{b_1}{2} + \frac{a}{z-z_0},\ b_1^2 - 4 b_0=0.
}
\end{array}
\right.
\label{eqRational}
\end{eqnarray}

\medskip
\noindent
{\bf Acknowledgement.} 
The authors would like to thank the referee for the valuable suggestions.


\vfill\eject
\end{document}